\documentclass[12pt]{article}
\usepackage{amssymb}
\usepackage[dvips]{graphicx}
\newcommand{\be}{\begin{equation}}
\newcommand{\ee}{\end{equation}}
\newcommand{\bea}{\begin{eqnarray}}
\newcommand{\eea}{\end{eqnarray}}
\newcommand{\ba}{\begin{array}}
\newcommand{\ea}{\end{array}}

\newcommand{\bc}{\begin{center}}
\newcommand{\ec}{\end{center}}
\newcommand{\ben}{\begin{enumerate}}
\newcommand{\een}{\end{enumerate}}
\newcommand{\bfi}{\begin{figure}}
\newcommand{\efi}{\end{figure}}

\newcommand{\bq}{\begin{quote}}
\newcommand{\eq}{\end{quote}}
\newcommand{\bqu}{\begin{quotation}}
\newcommand{\equ}{\end{quotation}}
\newenvironment{emphit}{\begin{itemize}}{\end{itemize}}
\newcommand{\bemp}{\begin{emphit}}
\newcommand{\eemp}{\end{emphit}}

\newcommand{\bt}{\begin{tabular}}
\newcommand{\et}{\end{tabular}}

\newtheorem{myth}{Theorem}[section]
\newtheorem{mylem}{Lemma}[section]

\newtheorem{mydef}{Definition}[section]
\newtheorem{myexam}{Example}[section]
\newtheorem{myprob}{Problem}[section]
\newtheorem{myquest}{Question}[section]
\newtheorem{myrem}{Remark}[section]

\begin{document}
\date{}
\title{Palindromic Permutations and Generalized Smarandache
Palindromic Permutations 
}
\author{T\`em\'it\'op\'e Gb\'ol\'ah\`an Ja\'iy\'e\d ol\'a\thanks{On Doctorate Programme at
the University of Agriculture Abeokuta, Nigeria.}\\Department of
Mathematics,\\
Obafemi Awolowo University, Ile Ife,
Nigeria.\\jaiyeolatemitope@yahoo.com, tjayeola@oauife.edu.ng}
\maketitle

\begin{abstract}
The idea of left(right) palindromic permutations(LPPs,RPPs) and
left(right) generalized Smarandache palindromic
permutations(LGSPPs,RGSPPs) are introduced in symmetric groups $S_n$
of degree $n$. It is shown that in $S_n$, there exist a LPP and a
RPP and they are unique(this fact is demonstrated using $S_2$ and
$S_3$). The dihedral group $D_n$ is shown to be generated by a RGSPP
and a LGSPP(this is observed to be true in $S_3$) but the geometric
interpretations of a RGSPP and a LGSPP are found not to be rotation
and reflection respectively. In $S_3$, each permutation is at least
a RGSPP or a LGSPP. There are $4$ RGSPPs and $4$ LGSPPs in $S_3$,
while $2$ permutations are both RGSPPs and LGSPPs. A permutation in
$S_n$ is shown to be a LPP or RPP(LGSPP or RGSPP) if and only if its
inverse is a LPP or RPP(LGSPP or RGSPP) respectively. Problems for
future studies are raised.
\end{abstract}

\section{Introduction}
\paragraph{}
According to Ashbacher and Neirynck \cite{gsp1}, an integer is said
to be a palindrome if it reads the same forwards and backwards. For
example, 12321 is a palindromic number. They also stated that it is
easy to prove that the density of the palindromes is zero in the set
of positive integers and they went ahead to answer the question on
the density of generalized Smarandache palindromes (GSPs) by showing
that the density of GSPs in the positive integers is approximately
0.11. Gregory \cite{gsp2}, Smarandache \cite{gsp3} and Ramsharan
\cite{gsp4} defined a generalized Smarandache palindrome (GSP) as
any integer or number of the form
\begin{displaymath}
a_1a_2a_3\cdots a_na_n\cdots a_3a_2a_1\qquad \textrm{or}\qquad
a_1a_2a_3\cdots a_{n-1}a_na_{n-1}\cdots a_3a_2a_1
\end{displaymath}
where all $a_1,a_2,a_3,\cdots a_n\in \mathbb{N}$ having one or more
digits. On the other hand, Hu \cite{gsp8} calls any integer or
number of this form a Smarandache generalized palindrome(SGP). His
naming will not be used here the first naming will be adopted.

Numbers of this form have also been considered by Khoshnevisan
\cite{gsp5}, \cite{gsp6} and \cite{gsp7}. For the sake of
clarification, it must be mentioned that the possibility of the
trivial case of enclosing the entire number is excluded. For
example, $12345$ can be written as $(12345)$. In this case, the
number is simply said to be a palindrome or a palindromic number as
it was mentioned earlier on. So, every number is a GSP. But this
possibility is eliminated by requiring that each number be split
into at least two segments if it is not a regular palindrome.
Trivially, since each regular palindrome is also a GSP and there are
GSPs that are not regular palindromes, there are more GSPs than
there are regular palindromes. As mentioned by Gregory \cite{gsp2},
very interesting GSPs are formed from smarandacheian sequences. For
an illustration he cited the smarandacheian sequence
\begin{displaymath}
11, 1221, 123321,\cdots , 123456789987654321,
1234567891010987654321,
\end{displaymath}
\begin{displaymath}
12345678910111110987654321, \cdots
\end{displaymath}
and observed that all terms are all GSPs. He also mentioned that it
has been proved that the GSP $1234567891010987654321$ is a prime and
concluded his work by possing the question of 'How many primes are
in the GSP sequence above?'.

Special mappings such as morphisms(homomorphisms, endomorphisms,
automorphisms, isomorphisms e.t.c) have been useful in the study of
the properties of most algebraic structures(e.g groupoids,
quasigroups, loops, semigroups, groups e.tc.). In this work, the
notion of palindromic permutations and generalized Smarandache
palindromic permutations are introduced and studied using the
symmetric group on the set $\mathbb{N}$ and this can now be viewed
as the study of some palindromes and generalized Smarandache
palindromes of numbers.

The idea of left(right) palindromic permutations(LPPs,RPPs) and
left(right) generalized Smarandache palindromic
permutations(LGSPPs,RGSPPs) are introduced in symmetric groups $S_n$
of degree $n$. It is shown that in $S_n$, there exist a LPP and a
RPP and they are unique. The dihedral group $D_n$ is shown to be
generated by a RGSPP and a LGSPP but the geometric interpretations
of a RGSPP and a LGSPP are found not to be rotation and reflection
respectively. In $S_3$, each permutation is at least a RGSPP or a
LGSPP. There are $4$ RGSPPs and $4$ LGSPPs in $S_3$, while $2$
permutations are both RGSPPs and LGSPPs. A permutation in $S_n$ is
shown to be a LPP or RPP(LGSPP or RGSPP) if and only if its inverse
is a LPP or RPP(LGSPP or RGSPP) respectively. Some of these results
are demonstrated with $S_2$ and $S_3$. Problems for future studies
are raised.

But before then, some definitions and basic results on symmetric
groups in classical group theory which shall be employed and used
are highlighted first.

\section{Preliminaries}
\begin{mydef}\label{1:1}
Let $X$ be a non-empty set. The group of all permutations of $X$
under composition of mappings is called the \textbf{symmetric group}
on $X$ and is denoted by $S_X$. A subgroup of $S_X$ is called a
permutation group on $X$.
\end{mydef}

It is easily seen that a bijection $X\simeq Y$ induces in a natural
way an isomorphism $S_X\cong S_Y$. If $|X|=n$, $S_X$ is denoted by
$S_n$ and called the {\it symmetric group of degree} $n$.

A permutation $\sigma\in S_n$ can be exhibited in the form

\begin{displaymath}
\left(\begin{array}{cccc}
1 & 2 & \cdots & n\\
\sigma (1) & \sigma (2) & \cdots & \sigma (n)
\end{array}\right),
\end{displaymath}
consisting of two rows of integers; the top row has integers
$1,2,\cdots ,n$ usually(but not necessarily) in their natural order,
and the bottom row has $\sigma (i)$ below $i$ for each $i=1,2,\cdots
,n$. This is called a two-row notation for a permutation. There is a
simpler, one-row notation for a special kind of permutation called
\textit{cycle}.

\begin{mydef}\label{1:2}
Let $\sigma\in S_n$. If there exists a list of distinct integers
$x_1,\cdots ,x_r\in \mathbb{N}$ such that
\begin{displaymath}
\begin{array}{lll}
\sigma (x_i)=x_{i+1},\qquad i=1,\cdots ,r-1,\\
\sigma (x_r)=x_1,\\
\sigma (x)=x~ \textrm{if $x\not \in \{x_1,\cdots ,x_r\}$},
\end{array}
\end{displaymath}
then $\sigma$ is called a cycle of length $r$ and denoted by
$(x_1\cdots x_r)$.
\end{mydef}
\begin{myrem}
A cycle of length $2$ is called a transposition. In other words , a
cycle $(x_1\cdots x_r)$ moves the integers $x_1,\cdots ,x_r$ one
step around a circle and leaves every other integer in $\mathbb{N}$.
If $\sigma (x)=x$, we say $\sigma$ does not move $x$. Trivially, any
cycle of length $1$ is the identity mapping $I$ or $e$. Note that
the one-row notation for a cycle does not indicate the degree $n$,
which has to be understood from the context.
\end{myrem}

\begin{mydef}\label{1:2.1}
Let $X$ be a set of points in space, so that the distance $d(x,y)$
between points $x$ and $y$ is given for all $x,y\in X$. A
permutation $\sigma$ of $X$ is called a \textbf{symmetry} of $X$ if
\begin{displaymath}
d(\sigma (x),\sigma (y))=d(x,y)~\forall~x,y\in X.
\end{displaymath}
Let $X$ be the set of points on the vertices of a regular polygon
which are labelled $\{1,2,\cdots,n\}$ i.e

The group of symmetries of a regular polygon $P_n$ of $n$ sides is
called the \textbf{dihedral group of  degree $n$} and denoted $D_n$.
\end{mydef}

\begin{myrem}
It must be noted that $D_n$ is a subgroup of $S_n$ i.e $D_n\le S_n$.
\end{myrem}

\begin{mydef}\label{1:3}
Let $S_n$ be a symmetric group of degree $n$. If $\sigma\in S_n$
such that
\begin{displaymath}
\sigma =\left(\begin{array}{cccc}
1 & 2 & \cdots & n\\
\sigma (1) & \sigma (2) & \cdots & \sigma (n)
\end{array}\right),
\end{displaymath}
then
\begin{enumerate}
\item the number $N_\lambda (\sigma)=12\cdots n\sigma (n)\cdots \sigma(1)$ is called the left palindromic value(LPV) of
$\sigma$.
\item the number $N_\rho (\sigma)=12\cdots n\sigma (1)\cdots \sigma(n)$ is called the right palindromic value(RPV) of
$\sigma$.
\end{enumerate}
\end{mydef}

\begin{mydef}\label{1:4}
Let $\sigma\in S_X$ such that
\begin{displaymath}
\sigma =\left(\begin{array}{cccc}
x_1 & x_2 & \cdots & x_n\\
\sigma (x_1) & \sigma (x_2) & \cdots & \sigma (x_n)
\end{array}\right).
\end{displaymath}
If $X=\mathbb{N}$, then
\begin{enumerate}
\item $\sigma$ is called a left palindromic permutation(LPP) if and
only if the number $N_\lambda (\sigma)$ is a palindrome.
\begin{displaymath}
PP_\lambda (S_X)=\{\sigma\in S_X~:~\textrm{$\sigma$ is a LPP}\}
\end{displaymath}
\item $\sigma$ is called a right palindromic permutation(RPP) if
and only if the number $N_\rho (\sigma)$ is a palindrome.
\begin{displaymath}
PP_\rho (S_X)=\{\sigma\in S_X~:~\textrm{$\sigma$ is a RPP}\}
\end{displaymath}
\item $\sigma$ is called a palindromic permutation(PP) if and
only if it is both a LPP and a RPP.
\begin{displaymath}
PP(S_X)=\{\sigma\in S_X~:~\textrm{$\sigma$ is a LPP and a RPP
}\}=PP_\lambda (S_X)\bigcap PP_\rho (S_X)
\end{displaymath}
\end{enumerate}
\end{mydef}

\begin{mydef}\label{1:5}
Let $\sigma\in S_X$ such that
\begin{displaymath}
\sigma =\left(\begin{array}{cccc}
x_1 & x_2 & \cdots & x_n\\
\sigma (x_1) & \sigma (x_2) & \cdots & \sigma (x_n)
\end{array}\right).
\end{displaymath}
If $X=\mathbb{N}$, then
\begin{enumerate}
\item $\sigma$ is called a left generalized Smarandache palindromic permutation(LGSPP) if and
only if the number $N_\lambda (\sigma)$ is a GSP.
\begin{displaymath}
GSPP_\lambda (S_X)=\{\sigma\in S_X~:~\textrm{$\sigma$ is a LGSPP}\}
\end{displaymath}
\item $\sigma$ is called a right generalized Smarandache palindromic permutation(RGSPP) if and
only if the number $N_\rho (\sigma)$ is a GSP.
\begin{displaymath}
GSPP_\rho (S_X)=\{\sigma\in S_X~:~\textrm{$\sigma$ is a RGSPP}\}
\end{displaymath}
\item $\sigma$ is called a generalized Smarandache palindromic permutation(GSPP) if and
only if it is both a LGSPP and a RGSPP.
\begin{displaymath}
GSPP(S_X)=\{\sigma\in S_X~:~\textrm{$\sigma$ is a LGSPP and a RGSPP
}\}=GSPP_\lambda (S_X)\bigcap GSPP_\rho (S_X)
\end{displaymath}
\end{enumerate}
\end{mydef}

\begin{myth}\label{1:5.1}(Cayley Theorem)

Every group is isomorphic to a permutation group.
\end{myth}

\begin{myth}\label{1:5.2}
The dihedral group $D_n$ is a group of order $2n$ generated by two
elements $\sigma,\tau$ satisfying $\sigma^n=e=\tau^2$ and
$\tau\sigma=\sigma^{n-1}\tau$, where
\begin{displaymath}
\sigma =\left(\begin{array}{cccc} 1 & 2 & \cdots & n
\end{array}\right)\qquad \textrm{and}\qquad
\tau =\left(\begin{array}{cccc}
1 & 2 & \cdots & n\\
1 & n & \cdots & 2
\end{array}\right).
\end{displaymath}
\end{myth}

\section{Main Results}
\begin{myth}\label{1:6}
In any symmetric group $S_n$ of degree $n$, there exists
\begin{enumerate}
\item a LPP and it is unique.
\item a RPP and it is unique.
\end{enumerate}
But there does not exist a PP.
\end{myth}
{\bf Proof}\\
Let $\sigma\in S_n$, then
\begin{displaymath}
\sigma =\left(\begin{array}{cccc}
x_1 & x_2 & \cdots & x_n\\
\sigma (x_1) & \sigma (x_2) & \cdots & \sigma (x_n)
\end{array}\right).
\end{displaymath}
\begin{enumerate}
\item When
\begin{displaymath}
\sigma (n)=n,\sigma (n-1)=n-1,\cdots,\sigma (2)=2,\sigma (1)=1
\end{displaymath}
then the number
\begin{displaymath}
N_\lambda (\sigma)=12\cdots n\sigma(n)\cdots \sigma(2)\sigma(1)
=12\cdots nn\cdots 21
\end{displaymath}
is a palindrome which implies $\sigma\in PP_\lambda (S_n)$. So,
there exists a LPP. The uniqueness is as follows. Observe that
\begin{displaymath}
\sigma =\left(\begin{array}{cccc}
1 & 2 & \cdots & n\\
1 & 2 & \cdots & n
\end{array}\right)=I.
\end{displaymath}
Since $S_n$ is a group for all $n\in \mathbb{N}$ and $I$ is the
identity element(mapping), then it must be unique.
\item When
\begin{displaymath}
\sigma (1)=n,\sigma (2)=n-1,\cdots,\sigma (n-1)=2,\sigma (n)=1
\end{displaymath}
then the number
\begin{displaymath}
N_\rho (\sigma)=12\cdots n\sigma(1)\cdots \sigma(n-1)\sigma(n)
=12\cdots nn\cdots 21
\end{displaymath}
is a palindrome which implies $\sigma\in PP_\rho (S_n)$. So, there
exists a RPP. The uniqueness is as follows. If there exist two of
such, say $\sigma_1$ and $\sigma_2$ in $S_n$, then
\begin{displaymath}
\sigma_1 =\left(\begin{array}{cccc}
1 & 2 & \cdots & n\\
\sigma_1 (1) & \sigma_1 (2) & \cdots & \sigma_1 (n)
\end{array}\right)\qquad \textrm{and}\qquad
\sigma_2 =\left(\begin{array}{cccc}
1 & 2 & \cdots & n\\
\sigma_2 (1) & \sigma_2 (2) & \cdots & \sigma_2 (n)
\end{array}\right)
\end{displaymath}
such that
\begin{displaymath}
N_\rho (\sigma_1)=12\cdots n\sigma_1(1)\cdots
\sigma_1(n-1)\sigma_1(n)
\end{displaymath}
and
\begin{displaymath}
N_\rho (\sigma_2)=12\cdots n\sigma_2(1)\cdots
\sigma_2(n-1)\sigma_2(n)
\end{displaymath}
are palindromes which implies
\begin{displaymath}
\sigma_1(1)=n,\sigma_1(2)=n-1,\cdots,\sigma_1(n-1)=2,\sigma_1(n)=1
\end{displaymath}
and
\begin{displaymath}
\sigma_2(1)=n,\sigma_2(2)=n-1,\cdots,\sigma_2(n-1)=2,\sigma_2(n)=1.
\end{displaymath}
So, $\sigma_1=\sigma_2$, thus $\sigma$ is unique.
\end{enumerate}
The proof of the last part is as follows. Let us assume by
contradiction that there exists a PP $\sigma\in S_n$. Then if
\begin{displaymath}
\sigma =\left(\begin{array}{cccc}
1 & 2 & \cdots & n\\
\sigma (1) & \sigma (2) & \cdots & \sigma (n)
\end{array}\right),
\end{displaymath}
\begin{displaymath}
N_\lambda (\sigma)=12\cdots n\sigma(n)\cdots \sigma(2)\sigma(1)
\end{displaymath}
and
\begin{displaymath}
N_\rho (\sigma)=12\cdots n\sigma(1)\cdots \sigma(n-1)\sigma(n)
\end{displaymath}
are palindromes. So that $\sigma\in S_n$ is a PP. Consequently,
\begin{displaymath}
n=\sigma (n)=1, n-1=\sigma (n-1)=2,\cdots, 1=\sigma (1)=n,
\end{displaymath}
so that $\sigma$ is not a bijection which means $\sigma\not \in
S_n$. This is a contradiction. Hence, no PP exist.

\begin{myexam}\label{1:7}
Let us consider the symmetric group $S_2$ of degree $2$. There are
two permutations of the set $\{1,2\}$ given by
\begin{displaymath}
I =\left(\begin{array}{cc}
1 & 2\\
1 & 2
\end{array}\right)\qquad \textrm{and}\qquad
\delta =\left(\begin{array}{cc}
1 & 2\\
2 & 1
\end{array}\right).
\end{displaymath}
\begin{displaymath}
N_\rho(I)=1212=(12)(12),N_\lambda (I)=1221~ \textrm{or}~ N_\lambda
(I)=1(22)1,
\end{displaymath}
\begin{displaymath}
N_\rho (\delta)=1221~  \textrm{or}~ N_\rho (\delta)=(12)(21)
~\textrm{and}~ N_\lambda (\delta)=1212=(12)(12).
\end{displaymath}
\end{myexam}
So, $I$ and $\delta$ are both RGSPPs and LGSPPs which implies $I$
and $\delta$ are GSPPs i.e $I,\delta\in GSPP_\rho (S_2)$ and
$I,\delta\in GSPP_\lambda (S_2)\Rightarrow I,\delta\in GSPP(S_2)$.
Therefore, $GSPP(S_2)=S_2$. Furthermore, it can be seen that the
result in Theorem~\ref{1:6} is true for $S_2$ because only $I$ is a
LPP and only $\delta$ is a RPP. There is definitely no PP as the
theorem says.

\begin{myexam}\label{1:8}
Let us consider the symmetric group $S_3$ of degree $3$. There are
six permutations of the set $\{1,2,3\}$ given by
\begin{displaymath}
e=I =\left(\begin{array}{ccc}
1 & 2 & 3\\
1 & 2 & 3
\end{array}\right),
\sigma_1 =\left(\begin{array}{ccc}
1 & 2 & 3\\
2 & 3 & 1
\end{array}\right),
\sigma_2 =\left(\begin{array}{ccc}
1 & 2 & 3\\
3 & 1 & 2
\end{array}\right),
\end{displaymath}
\begin{displaymath}
\tau_1 =\left(\begin{array}{ccc}
1 & 2 & 3\\
1 & 3 & 2
\end{array}\right),
\tau_2 =\left(\begin{array}{ccc}
1 & 2 & 3\\
3 & 2 & 1
\end{array}\right)\qquad \textrm{and}\qquad
\tau_3 =\left(\begin{array}{ccc}
1 & 2 & 3\\
2 & 1 & 3
\end{array}\right).
\end{displaymath}
As claimed in Theorem~\ref{1:6}, the unique LPP in $S_3$ is $I$
while the unique RPP in $S_3$ is $\tau_2$. There is no PP as the
theorem says.
\end{myexam}

\begin{mylem}\label{1:9}
In $S_3$, the following are true.
\begin{enumerate}
\item At least $\sigma\in GSPP_\rho (S_3)$ or $\sigma\in
GSPP_\lambda (S_3)~\forall~\sigma\in S_3$.
\item $|GSPP_\rho (S_3)|=4$, $|GSPP_\lambda (S_3)|=4$ and $|GSPP(S_3)|=2$.
\end{enumerate}
\end{mylem}
{\bf Proof}\\
Observe the following :
\begin{displaymath}
N_\lambda (I)=123321,~N_\rho (I)=123123=(123)(123).
\end{displaymath}
\begin{displaymath}
N_\lambda (\sigma_1)=123132,~N_\rho (\sigma_1)=123231=1(23)(23)1.
\end{displaymath}
\begin{displaymath}
N_\lambda (\sigma_2)=123213,~N_\rho (\sigma_2)=123312=(12)(33)(12).
\end{displaymath}
\begin{displaymath}
N_\lambda (\tau_1)=123231=1(23)(23)1,~N_\rho (\tau_1)=123132.
\end{displaymath}
\begin{displaymath}
N_\lambda (\tau_2)=123123=(123)(123),~N_\rho (\tau_2)=123321=123321.
\end{displaymath}
\begin{displaymath}
N_\lambda (\tau_3)=123312=(12)(33)(12),~N_\rho (\tau_3)=123213.
\end{displaymath}
So, $GSPP_\lambda (S_3)=\{I,\tau_1,\tau_2,\tau_3\}$ and $GSPP_\rho
(S_3)=\{I,\sigma_1,\sigma_2,\tau_2\}$. Thus, 1. is true. Therefore,
$|GSPP_\rho (S_3)|=4$, $|GSPP_\lambda (S_3)|=4$ and
$|GSPP(S_3)|=|GSPP_\rho (S_3)\bigcap GSPP_\lambda (S_3)|=2$. So, 2.
is true.

\begin{mylem}\label{1:10}
$S_3$ is generated by a RGSPP and a LGSPP.
\end{mylem}
{\bf Proof}\\
Recall from Example~\ref{1:8} that
\begin{displaymath}
S_3=\{I=e,\sigma_1,\sigma_2,\tau_1,\tau_2,\tau_3\}.
\end{displaymath}
If $\sigma =\sigma_1$ and $\tau =\tau_1$, then it is easy to verify
that
\begin{displaymath}
\sigma^2=\sigma_2,~\sigma^3=e,~\tau^2=e,~\sigma\tau
=\tau_3,~\sigma^2\tau
=\tau_2=\tau\sigma~\textrm{hence,}
\end{displaymath}
\begin{displaymath}
S_3=\{e,\sigma,\sigma^2,\tau,\sigma\tau,\sigma^2\tau_3\}\Rightarrow
S_3=\big<\sigma ,\tau\big>.
\end{displaymath}
From the proof Lemma~\ref{1:9}, $\sigma$ is a RGSPP and $\tau$ is a
LGSPP. This justifies the claim.

\begin{myrem}
In Lemma~\ref{1:10}, $S_3$ is generated by a RGSPP and a LGSPP.
Could this statement be true for all $S_n$ of degree $n$? Or could
it be true for some subgroups of $S_n$? Also, it is interesting to
know the geometric meaning of a RGSPP and a LGSPP. So two questions
are possed and the two are answered.
\end{myrem}

\begin{myquest}\label{1:11}
\begin{enumerate}
\item Is the symmetric group $S_n$ of degree $n$ generated by a RGSPP
and a LGSPP? If not, what permutation group(s) is generated by a
RGSPP and a LGSPP?
\item Are the geometric interpretations of a RGSPP and a LGSPP rotation and reflection respectively?
\end{enumerate}
\end{myquest}

\begin{myth}\label{1:12}

The dihedral group $D_n$ is generated by a RGSPP and a LGSPP i.e
$D_n=\big<\sigma,\tau\big>$ where $\sigma\in GSPP_\rho (S_n)$ and
$\tau\in GSPP_\lambda (S_n)$.
\end{myth}
{\bf Proof}\\
Recall from Theorem\ref{1:5.2} that the dihedral group
$D_n=\big<\sigma,\tau\big>$ where
\begin{displaymath}
\sigma =\left(\begin{array}{cccc} 1 & 2 & \cdots & n
\end{array}\right)=\left(\begin{array}{cccc} 1 & 2 & \cdots & n\\
2 & 3 & \cdots & 1
\end{array}\right)\qquad \textrm{and}\qquad
\tau =\left(\begin{array}{cccc}
1 & 2 & \cdots & n\\
1 & n & \cdots & 2
\end{array}\right).
\end{displaymath}
Observe that
\begin{displaymath}
N_\rho (\sigma )=123\cdots n23\cdots n1=1(23\cdots n)(23\cdots
n)1,~N_\lambda (\sigma)=123\cdots n1n\cdots 32.
\end{displaymath}
\begin{displaymath}
N_\rho (\tau )=12\cdots n1n\cdots 2,~N_\lambda (\tau)=12\cdots
n2\cdots n1=1(2\cdots n)(2\cdots n)1.
\end{displaymath}
So, $\sigma\in GSPP_\rho (S_n)$ and $\tau\in GSPP_\lambda (S_n)$.
Therefore, the dihedral group $D_n$ is generated by a RGSPP and a
LGSPP.

\begin{myrem}
In Lemma~\ref{1:10}, it was shown that $S_3$ is generated by a RGSPP
and a LGSPP. Considering Theorem~\ref{1:12} when $n=3$, it can be
deduced that $D_3$ will be generated by a RGSPP and a LGSPP. Recall
that $|D_3|=2\times 3=6$, so $S_3=D_3$. Thus Theorem~\ref{1:12}
generalizes Lemma~\ref{1:10}.
\end{myrem}

\paragraph{Rotations and Reflections}
Geometrically, in Theorem~\ref{1:12}, $\sigma$ is a rotation of the
regular polygon $P_n$ through an angle $\frac{2\pi}{n}$ in its own
plane, and $\tau$ is a reflection (or a turning over) in the
diameter through the vertex $1$. It looks like a RGSPP and a LGSPP
are formed by rotation and reflection respectively. But there is a
contradiction in $S_4$ which can be traced from a subgroup of $S_4$
particularly the Klein four-group. The Klein four-group is the group
of symmetries of a four sided non-regular polygon(rectangle). The
elements are:
\begin{displaymath}
e=I =\left(\begin{array}{cccc}
1 & 2 & 3 & 4\\
1 & 2 & 3 & 4
\end{array}\right),
\delta_1 =\left(\begin{array}{cccc}
1 & 2 & 3 & 4\\
3 & 4 & 1 & 2
\end{array}\right),
\delta_2 =\left(\begin{array}{cccc}
1 & 2 & 3 & 4\\
2 & 1 & 4 & 3
\end{array}\right)
\end{displaymath}
\begin{displaymath}
\qquad\textrm{and}\qquad \delta_3 =\left(\begin{array}{cccc}
1 & 2 & 3 & 4\\
4 & 3 & 2 & 1
\end{array}\right).
\end{displaymath}
Observe the following:
\begin{displaymath}
N_\rho (\delta_1)=12343412=(12)(34)(34)(12),~N_\lambda
(\delta_1)=12342143.
\end{displaymath}
\begin{displaymath}
N_\rho (\delta_2)=12342143=12342143,~N_\lambda
(\delta_2)=12343412=(12)(34)(34)(12).
\end{displaymath}
\begin{displaymath}
N_\rho (\delta_3)=12344321=123(44)321,~N_\lambda
(\delta_3)=12341234=(1234)(1234).
\end{displaymath}
So, $\delta_1$ is a RGSPP while $\delta_2$ is a LGSPP and $\delta_3$
is a GSPP. Geometrically, $\delta_1$ is a rotation through an angle
of $\pi$ while $\delta_2$ and $\delta_3$ are reflections in the axes
of symmetry parallel to the sides. Thus $\delta_3$ which is a GSPP
is both a reflection and a rotation, which is impossible. Therefore,
the geometric meaning of a RGSPP and a LGSPP are not rotation and
reflection respectively. It is difficult to really ascertain the
geometric meaning of a RGSPP and a LGSPP if at all it exist.

\paragraph{}
How beautiful will it be if $GSPP_\rho (S_n)$, $PP_\rho (S_n)$,
$GSPP_\lambda (S_n)$, $PP_\lambda (S_n)$, $GSPP(S_n)$ and $PP(S_n)$
form algebraic structures under the operation of map composition.

\begin{myth}\label{1:13}
Let $S_n$ be a symmetric group of degree $n$. If $\sigma\in S_n$,
then
\begin{enumerate}
\item $\sigma\in PP_\lambda (S_n)\Leftrightarrow \sigma^{-1}\in PP_\lambda
(S_n)$.
\item
$\sigma\in PP_\rho (S_n)\Leftrightarrow \sigma^{-1}\in PP_\rho
(S_n)$.
\item $I\in PP_\lambda (S_n)$.
\end{enumerate}
\end{myth}
{\bf Proof}\\
\begin{enumerate}
\item $\sigma\in PP_\lambda (S_n)$ implies
\begin{displaymath}
N_\lambda (\sigma)=12\cdots n\sigma(n)\cdots \sigma(2)\sigma(1)
\end{displaymath}
is a palindrome. Consequently,
\begin{displaymath}
\sigma (n)=n,\sigma (n-1)=n-1,\cdots,\sigma (2)=2,\sigma (1)=1.
\end{displaymath}
So,
\begin{displaymath}
N_\lambda (\sigma^{-1})=\sigma(1)\sigma(2)\cdots \sigma(n)n\cdots
21=12\cdots nn\cdots 21\Rightarrow \sigma^{-1}\in PP_\lambda (S_n).
\end{displaymath}
The converse is similarly proved by carrying out the reverse of the
procedure above.
\item $\sigma\in PP_\rho (S_n)$ implies
\begin{displaymath}
N_\rho (\sigma)=12\cdots n\sigma(1)\cdots \sigma(n-1)\sigma(n)
\end{displaymath}
is a palindrome. Consequently,
\begin{displaymath}
\sigma (1)=n,\sigma (2)=n-1,\cdots,\sigma (n-1)=2,\sigma (n)=1.
\end{displaymath}
So,
\begin{displaymath}
N_\rho (\sigma^{-1})=\sigma(1)\cdots \sigma(n-1)\sigma(n)12\cdots
n=n\cdots 2112\cdots n\Rightarrow \sigma^{-1}\in PP_\rho (S_n).
\end{displaymath}
The converse is similarly proved by carrying out the reverse of the
procedure above.
\item
\begin{displaymath}
I =\left(\begin{array}{cccc}
1 & 2 & \cdots & n\\
1 & 2 & \cdots & n
\end{array}\right).
\end{displaymath}
\begin{displaymath}
N_\lambda (I)=12\cdots nn\cdots 21\Rightarrow I\in PP_\lambda (S_n).
\end{displaymath}
\end{enumerate}

\begin{myth}\label{1:14}
Let $S_n$ be a symmetric group of degree $n$. If $\sigma\in S_n$,
then
\begin{enumerate}
\item $\sigma\in GSPP_\lambda (S_n)\Leftrightarrow \sigma^{-1}\in GSPP_\lambda
(S_n)$.
\item $\sigma\in GSPP_\rho (S_n)\Leftrightarrow \sigma^{-1}\in GSPP_\rho
(S_n)$.
\item $I\in GSPP (S_n)$.
\end{enumerate}
\end{myth}
{\bf Proof}\\
If $\sigma\in S_n$, then
\begin{displaymath}
\sigma =\left(\begin{array}{cccc}
1 & 2 & \cdots & n\\
\sigma (1) & \sigma (2) & \cdots & \sigma (n)
\end{array}\right).
\end{displaymath}
So,
\begin{displaymath}
N_\lambda (\sigma)=12\cdots n\sigma(n)\cdots \sigma(2)\sigma(1)
\end{displaymath}
and
\begin{displaymath}
N_\rho (\sigma)=12\cdots n\sigma(1)\cdots \sigma(n-1)\sigma(n)
\end{displaymath} are numbers with even number of digits whether $n$
is an even or odd number. Thus, $N_\rho (\sigma)$ and $N_\lambda
(\sigma)$ are GSPs defined by
\begin{displaymath}
a_1a_2a_3\cdots a_na_n\cdots a_3a_2a_1
\end{displaymath}
and not
\begin{displaymath}
a_1a_2a_3\cdots a_{n-1}a_na_{n-1}\cdots a_3a_2a_1
\end{displaymath}
where all $a_1,a_2,a_3,\cdots a_n\in \mathbb{N}$ having one or more
digits because the first has even number of digits(or grouped
digits) while the second has odd number of digits(or grouped
digits). The following grouping notations will be used:
\begin{displaymath}
(a_i)_{i=1}^n=a_1a_2a_3\cdots a_n\qquad \textrm{and}\qquad
[a_i]_{i=1}^n=a_na_{n-1}a_{n-2}\cdots a_3a_2a_1.
\end{displaymath}
Let $\sigma\in S_n$ such that
\begin{displaymath}
\sigma =\left(\begin{array}{cccc}
x_1 & x_2 & \cdots & x_n\\
\sigma (x_1) & \sigma (x_2) & \cdots & \sigma (x_n)
\end{array}\right)
\end{displaymath}
where $x_i\in \mathbb{N}~\forall~i\in \mathbb{N}$.
\begin{enumerate}
\item So, $\sigma\in GSPP_\lambda (S_n)$ implies
\begin{displaymath}
N_\lambda (\sigma
)=(x_{i_1})_{i_1=1}^{n_1}(x_{i_2})_{i_2=(n_1+1)}^{n_2}(x_{i_3})_{i_3=(n_2+
1)}^{n_3}\cdots
(x_{i_{n-1}})_{i_{n-1}=(n_{n-2}+1)}^{n_{n-1}}(x_{i_n})_{i_n=(n_{n-1}+1)}^{n_n}
\end{displaymath}
\begin{displaymath}\leftrightarrow\Leftrightarrow
[\sigma(x_{i_n})]_{i_n=(n_{n-1}+1)}^{n_n}
[\sigma(x_{i_{n-1}})]_{i_{n-1}=(n_{n-2}+1)}^{n_{n-1}}\cdots
[\sigma(x_{i_3})]_{i_3=(n_2+
1)}^{n_3}[\sigma(x_{i_2})]_{i_2=(n_1+1)}^{n_2}[\sigma(x_{i_1})]_{i_1=1}^{n_1}
\end{displaymath}
is a GSP where $x_{i_j}\in  \mathbb{N}~\forall~i_j\in \mathbb{N}$,
$j\in \mathbb{N}$ and $n_n=n$. The interval of integers [1,n] is
partitioned into
\begin{displaymath}
[1,n]=[1,n_1]\cup [n_1+1,n_2]\cup\cdots\cup [n_{n-2}+1,n_{n-1}]\cup
[n_{n-1},n_n].
\end{displaymath}
The length of each grouping $(\cdot )_{i_j}^{n_j}$ or $[\cdot
]_{i_j}^{n_j}$ is determined by the corresponding interval of
integers $[n_i+1,n_{i+1}]$ and it is a matter of choice in other to
make the number $N_\lambda (\sigma )$ a GSP.

Now that $N_\lambda (\sigma )$ is a GSP, the following are true:
\begin{displaymath}
(x_{i_n})_{i_n=(n_{n-1}+1)}^{n_n}=[\sigma(x_{i_n})]_{i_n=(n_{n-1}+1)}^{n_n}\qquad\Leftrightarrow\qquad
[x_{i_n}]_{i_n=(n_{n-1}+1)}^{n_n}=(\sigma(x_{i_n}))_{i_n=(n_{n-1}+1)}^{n_n}
\end{displaymath}
\begin{displaymath}
(x_{i_{n-1}})_{i_{n-1}=(n_{n-2}+1)}^{n_{n-1}}=[\sigma(x_{i_{n-1}})]_{i_{n-1}=(n_{n-2}+1)}^{n_{n-1}}~\Leftrightarrow~
[x_{i_{n-1}}]_{i_{n-1}=(n_{n-2}+1)}^{n_{n-1}}=(\sigma(x_{i_{n-1}}))_{i_{n-1}=(n_{n-2}+1)}^{n_{n-1}}
\end{displaymath}
\begin{displaymath}
\vdots\qquad\qquad\qquad\qquad\qquad\qquad\vdots\qquad\qquad\qquad\qquad\qquad\qquad\vdots
\end{displaymath}
\begin{displaymath}
(x_{i_2})_{i_2=(n_1+1)}^{n_2}=[\sigma(x_{i_2})]_{i_2=(n_1+1)}^{n_2}\qquad\Leftrightarrow\qquad
[x_{i_2}]_{i_2=(n_1+1)}^{n_2}=(\sigma(x_{i_2}))_{i_2=(n_1+1)}^{n_2}
\end{displaymath}
\begin{displaymath}
(x_{i_1})_{i_1=1}^{n_1}=[\sigma(x_{i_1})]_{i_1=1}^{n_1}\qquad\Leftrightarrow\qquad
[x_{i_1}]_{i_1=1}^{n_1}=(\sigma(x_{i_1}))_{i_1=1}^{n_1}
\end{displaymath}
Therefore, since
\begin{displaymath}
\sigma =\left(\begin{array}{ccccccccccc} x_1 & \cdots & x_{i_1} &
\cdots & x_{n_1} & \cdots & x_{n_{n-1}+1} & \cdots & x_{j_k} &
\cdots &
x_{n_n}  \\
\sigma(x_1) & \cdots & \sigma(x_{i_1}) & \cdots & \sigma(x_{n_1}) &
\cdots & \sigma(x_{n_{n-1}+1}) & \cdots & \sigma(x_{j_k}) & \cdots &
\sigma(x_{n_n})
\end{array}\right),
\end{displaymath}
then
\begin{displaymath}
\sigma^{-1} =\left(\begin{array}{ccccccccccc} \sigma(x_1) & \cdots &
\sigma(x_{i_1}) & \cdots & \sigma(x_{n_1}) & \cdots &
\sigma(x_{n_{n-1}+1}) & \cdots & \sigma(x_{j_k}) & \cdots &
\sigma(x_{n_n})\\
x_1 & \cdots & x_{i_1} & \cdots & x_{n_1} & \cdots & x_{n_{n-1}+1} &
\cdots & x_{j_k} & \cdots &
x_{n_n}  \\
\end{array}\right),
\end{displaymath}
so
\begin{displaymath}
N_\lambda (\sigma^{-1}
)=(\sigma(x_{i_1}))_{i_1=1}^{n_1}(\sigma(x_{i_2}))_{i_2=(n_1+1)}^{n_2}(\sigma(x_{i_3}))_{i_3=(n_2+
1)}^{n_3}\cdots
(\sigma(x_{i_{n-1}}))_{i_{n-1}=(n_{n-2}+1)}^{n_{n-1}}
\end{displaymath}
\begin{displaymath}
(\sigma(x_{i_n}))_{i_n=(n_{n-1}+1)}^{n_n}[x_{i_n}]_{i_n=(n_{n-1}+1)}^{n_n}[x_{i_{n-1}}]_{i_{n-1}=(n_{n-2}+1)}^{n_{n-1}}\cdots
[x_{i_3}]_{i_3=(n_2+
1)}^{n_3}[x_{i_2}]_{i_2=(n_1+1)}^{n_2}[x_{i_1}]_{i_1=1}^{n_1}
\end{displaymath}
is a GSP hence, $\sigma^{-1}\in GSPP_\lambda (S_n)$.

The converse can be proved in a similar way since
$(\sigma^{-1})^{-1}=\sigma$.

\item Also, $\sigma\in GSPP_\rho (S_n)$ implies
\begin{displaymath}
N_\rho (\sigma
)=(x_{i_1})_{i_1=1}^{n_1}(x_{i_2})_{i_2=(n_1+1)}^{n_2}(x_{i_3})_{i_3=(n_2+
1)}^{n_3}\cdots
(x_{i_{n-1}})_{i_{n-1}=(n_{n-2}+1)}^{n_{n-1}}(x_{i_n})_{i_n=(n_{n-1}+1)}^{n_n}
\end{displaymath}
\begin{displaymath}
(\sigma(x_{i_1}))_{i_1=1}^{n_1}(\sigma(x_{i_2}))_{i_2=(n_1+1)}^{n_2}(\sigma(x_{i_3}))_{i_3=(n_2+
1)}^{n_3}\cdots
(\sigma(x_{i_{n-1}}))_{i_{n-1}=(n_{n-2}+1)}^{n_{n-1}}(\sigma(x_{i_n}))_{i_n=(n_{n-1}+1)}^{n_n}
\end{displaymath}
is a GSP where $x_{i_j}\in  \mathbb{N}~\forall~i_j\in \mathbb{N}$,
$j\in \mathbb{N}$ and $n_n=n$. The interval of integers [1,n] is
partitioned into
\begin{displaymath}
[1,n]=[1,n_1]\cup [n_1+1,n_2]\cup\cdots\cup [n_{n-2}+1,n_{n-1}]\cup
[n_{n-1},n_n].
\end{displaymath}
The length of each grouping $(\cdot )_{i_j}^{n_j}$ is determined by
the corresponding interval of integers $[n_i+1,n_{i+1}]$ and it is a
matter of choice in other to make the number $N_\rho (\sigma )$ a
GSP.

Now that $N_\rho (\sigma )$ is a GSP, the following are true:
\begin{displaymath}
(x_{i_n})_{i_n=(n_{n-1}+1)}^{n_n}=(\sigma(x_{i_1}))_{i_1=1}^{n_1}
\end{displaymath}
\begin{displaymath}
(x_{i_{n-1}})_{i_{n-1}=(n_{n-2}+1)}^{n_{n-1}}=(\sigma(x_{i_2}))_{i_2=(n_1+1)}^{n_2}
\end{displaymath}
\begin{displaymath}
\vdots\qquad\qquad\qquad\vdots\qquad\qquad\qquad\vdots
\end{displaymath}
\begin{displaymath}
(x_{i_2})_{i_2=(n_1+1)}^{n_2}=(\sigma(x_{i_{n-1}}))_{i_{n-1}=(n_{n-2}+1)}^{n_{n-1}}
\end{displaymath}
\begin{displaymath}
(x_{i_1})_{i_1=1}^{n_1}=(\sigma(x_{i_n}))_{i_n=(n_{n-1}+1)}^{n_n}
\end{displaymath}
Therefore, since
\begin{displaymath}
\sigma =\left(\begin{array}{ccccccccccc} x_1 & \cdots & x_{i_1} &
\cdots & x_{n_1} & \cdots & x_{n_{n-1}+1} & \cdots & x_{j_k} &
\cdots &
x_{n_n}  \\
\sigma(x_1) & \cdots & \sigma(x_{i_1}) & \cdots & \sigma(x_{n_1}) &
\cdots & \sigma(x_{n_{n-1}+1}) & \cdots & \sigma(x_{j_k}) & \cdots &
\sigma(x_{n_n})
\end{array}\right),
\end{displaymath}
then
\begin{displaymath}
\sigma^{-1} =\left(\begin{array}{ccccccccccc} \sigma(x_1) & \cdots &
\sigma(x_{i_1}) & \cdots & \sigma(x_{n_1}) & \cdots &
\sigma(x_{n_{n-1}+1}) & \cdots & \sigma(x_{j_k}) & \cdots &
\sigma(x_{n_n})\\
x_1 & \cdots & x_{i_1} & \cdots & x_{n_1} & \cdots & x_{n_{n-1}+1} &
\cdots & x_{j_k} & \cdots &
x_{n_n}  \\
\end{array}\right),
\end{displaymath}
so
\begin{displaymath}
N_\rho
(\sigma^{-1})=(\sigma(x_{i_1}))_{i_1=1}^{n_1}(\sigma(x_{i_2}))_{i_2=(n_1+1)}^{n_2}(\sigma(x_{i_3}))_{i_3=(n_2+
1)}^{n_3}\cdots
(\sigma(x_{i_{n-1}}))_{i_{n-1}=(n_{n-2}+1)}^{n_{n-1}}
\end{displaymath}
\begin{displaymath}
(\sigma(x_{i_n}))_{i_n=(n_{n-1}+1)}^{n_n}(x_{i_1})_{i_1=1}^{n_1}(x_{i_2})_{i_2=(n_1+1)}^{n_2}(x_{i_3})_{i_3=(n_2+
1)}^{n_3}\cdots
(x_{i_{n-1}})_{i_{n-1}=(n_{n-2}+1)}^{n_{n-1}}(x_{i_n})_{i_n=(n_{n-1}+1)}^{n_n}
\end{displaymath}

is a GSP hence, $\sigma^{-1}\in GSPP_\rho (S_n)$.

The converse can be proved in a similar way since
$(\sigma^{-1})^{-1}=\sigma$.
\item
\begin{displaymath}
I =\left(\begin{array}{cccc}
1 & 2 & \cdots & n\\
1 & 2 & \cdots & n
\end{array}\right).
\end{displaymath}
\begin{displaymath}
N_\lambda (I)=12\cdots nn\cdots 21=12\cdots (nn)\cdots 21\Rightarrow
I\in GSPP_\lambda (S_n)~ \textrm{and}
\end{displaymath}
\begin{displaymath}
N_\rho (I)=(12\cdots n)(12\cdots n)\Rightarrow I\in GSPP_\rho (S_n)
\end{displaymath}
thus, $I\in GSPP (S_n)$.
\end{enumerate}

\section{Conclusion and Future studies}
By Theorem~\ref{1:6}, it is certainly true in every symmetric group
$S_n$ of degree $n$ there exist at least a RGSPP and a
LGSPP(although they are actually RPP and LPP). Following
Example~\ref{1:7}, there are $2$ RGSPPs, $2$ LGSPPs and $2$ GSPPs in
$S_2$ while from Lemma~\ref{1:9}, there are $4$ RGSPPs, $4$ LGSPPs
and $2$ GSPPs in $S_3$. Also, it can be observed that
\begin{displaymath}
|GSPP_\rho (S_2)|+|GSPP_\lambda
(S_2)|-|GSPP(S_2)|=2!=|S_2|~\textrm{and}
\end{displaymath}
\begin{displaymath}
|GSPP_\rho (S_3)|+|GSPP_\lambda (S_3)|-|GSPP(S_3)|=3!=|S_3|.
\end{displaymath}
The following problems are open for further studies.
\begin{myprob}
\begin{enumerate}
\item How many RGSPPs, LGSPPs and GSPPs are in $S_n$?
\item Does there exist functions $f_1,f_2,f_3~:~\mathbb{N}\to \mathbb{N}$ such that
$|GSPP_\rho (S_n)|=f_1(n)$, $|GSPP_\lambda (S_n)|=f_2(n)$ and
$|GSPP(S_n)|=f_3(n)$?
\item In general, does the formula
\begin{displaymath}
|GSPP_\rho (S_n)|+|GSPP_\lambda (S_n)|-|GSPP(S_n)|=n!=|S_n|?
\end{displaymath}
hold. If not, for what other $n>3$ is it true?
\end{enumerate}
\end{myprob}
The GAP package or any other appropriate mathematical package could
be helpful in investigating the solutions to them.

If the first question is answered, then the number of palindromes
that can be formed from the set $\{1,2,\cdots n\}$ can be known
since in the elements of $S_n$, the bottom row gives all possible
permutation of the integers $1,2,\cdots n$.

The Cayley Theorem(Theorem\ref{1:5.1}) can also be used to make a
further study on generalized Smarandache palindromic permutations.
In this work, $\mathbb{N}$ was the focus and it does not contain the
integer zero. This weakness can be strengthened by considering the
set $\mathbb{Z}_n=\{0,1,2,\cdots n-1\}~\forall~n\in \mathbb{N}$.
Recall that $(\mathbb{Z}_n,+)$ is a group and so by
Theorem~\ref{1:5.1} $(\mathbb{Z}_n,+)$ is isomorphic to a
permutation group particularly, one can consider a subgroup of the
symmetric group $S_{\mathbb{Z}_n}$.


\begin{thebibliography}{99}
\bibitem{gsp1} C. Ashbacher, L. Neirynck, {\it The density of generalized
Smarandache palindromes}, http://www.gallup.unm.edu/$\sim$
smarandache/GeneralizedPalindromes.htm.
\bibitem{gsp2} G. Gregory, {\it Generalized Smarandache palindromes}, http://www.gallup.unm.edu/$\sim$ smarandache/GSP.htm.
\bibitem{gsp8} C. Hu (2000), {\it On Smarandache generalized palindrome}, Second International Conference on Smarandache Type Notions In Mathematics and Quantum
Physics, University of Craiova, Craiova Romania,
atlas-conferences.com/c/a/f/t/25.htm.
\bibitem{gsp5} M. Khoshnevisan (2003), Manuscript.
\bibitem{gsp6} M. Khoshnevisan (2003), {\it Generalized Smarandache palindrome}, Mathematics
Magazine, Aurora, Canada.
\bibitem{gsp7} M. Khoshnevisan (2003), {\it Proposed problem 1062},
The PME Journal, USA, Vol. 11, No. 9, p. 501.
\bibitem{gsp4} K. Ramsharan (2003), Manusript, www.research.att.com/$\sim$njas/sequences/A082461.
\bibitem{gsp3} F. Smarandache (2006), {\it Sequences of Numbers Involved in
unsolved problems}, http://www.gallup.unm.edu/$\sim$
smarandache/eBooks-otherformats.htm, 140pp.
\end{thebibliography}
\end{document}